\documentclass[11pt,a4paper,twoside]{amsart}

\usepackage{amsfonts, indentfirst}
\usepackage{ifthen,amssymb}
\usepackage[all]{xypic}

\newtheorem{prop}{Proposition}[section]
\newtheorem{theorem}[prop]{Theorem}
\newtheorem{lemma}[prop]{Lemma}
\newtheorem{example}{Example}
\newtheorem{defi}[prop]{Definition}

\newtheorem{remark}[prop]{Remark}

\newcommand{\cqd}{\hfill$\Box$}
\newcommand{\mor}[0]{\operatorname{Mor}}
\newcommand{\quot}[0]{\operatorname{Quot}}
\renewcommand{\deg}[0]{{\it deg}}
\renewcommand{\dim}[0]{\operatorname{dim}}
\newcommand{\codim}[0]{\operatorname{codim}}
\newcommand{\im}[0]{\operatorname{Im}}
\newcommand{\ch}[0]{\operatorname{ch}}
\newcommand{\Td}[0]{\operatorname{Td}}
\newcommand{\rank}[0]{\operatorname{rank}}

\author[Cristina Mart{\'\i}nez]{Cristina Mart{\'\i}nez}
\author[Alberto Besana]{Alberto Besana}

\address{
Departament de Matematiques, Edifici C, Facultad de Ciencies, Universitat Aut\`onoma de Barcelona, 08193 Bellaterra, Barcelona}

\email{abesana@amaris.cat}
 \email{abesana@mat.uab.cat}
\address{ IIIA-Institut d'investigaci\'on en Intel$\cdot$ig\'encia Artificial, CSIC- Consejo Superior de Investigaciones Cient\'{\i}ficas, Campus UAB s/n, 08193 Bellaterra}
\email{cmartine@mat.uab.cat}
 \email{cristina@iiia.csic.es}
\title[Some remarks on cyclic Galois coverings of the projective line over finite fields]
{Some remarks on cyclic Galois coverings of the projective line over finite fields}
\begin{document}
\maketitle

\begin{abstract}
We study finite cyclic Galois extensions of the rational function field of  the projective line $\mathbb{P}^{1}(\mathbb{F}_{q})$ over a finite field $\mathbb{F}_{q}$ of $q$ elements defined by considering quotient curves by finite subgroups of the projective linear group $PGL(2,q)$,  and we enumerate them expressing the count in terms of Stirling numbers.

\end{abstract}

\section{Introduction}
Let $\mathbb{F}_{q}$ be the finite field with $q=p^{n}$ elements where $p$ is a prime number and $n\geq 1$ an integer.  Any other field $F$ of characteristic $p$ contains a copy of $\mathbb{F}_{p}$. We denote respectively by $\mathbb{A}^{n}(\mathbb{F}_{q})$ and $\mathbb{P}^{n}(\mathbb{F}_{q})$ the affine space and the projective space over $\mathbb{F}_{q}$. 
%Let denote by $\mathbb{F}_{p}$ the Galois field of $p$ elements. Any other field $F$ of characteristic $p$ contains a copy of $\mathbb{F}_{p}$. Any $V=\mathbb{F}_{p^{n}}$ field extension of $\mathbb{F}_{p}$ is a $\mathbb{F}_{p}$ vector space and a $n-1$ dimensional projective space  $\mathbb{P}^{n-1}(\mathbb{F}_{p})$. %We denote by $PG^{r}(n,p)$, the set of $r-The multiplication map
% $m_{y}: \mathbb{F}^{*}_{p}\rightarrow \mathbb{F}^{*}_{p}$ , mapping $x\mapsto y\,x$ is $\mathbb{F}_{p^{2}}$-linear and induces an automorphism of the projective plane $\mathbb{P}^{2}(\mathbb{F}_{p^{2}})$ of order $p^{2}$. 
Two points $x=(x_{0},x_{1},x_{2})$ and $y=(y_{0},y_{1},y_{2}) \in \mathbb{F}_{p^{3}}$ are projectively equivalent if and only if $y_{i}=t\,x_{i}$ for $i\in \{0,1,2\}$,  for some $t\in \mathbb{F}^{*}_{p}$.
Therefore the number of points in $\mathbb{P}^{2}(\mathbb{F}_{p})$ is $\frac{p^{3}-1}{p-1}=p^{2}+p+1$ and dually there are $p^{2}+p+1$ lines in $\mathbb{P}^{2}(\mathbb{F}_{p})$. %The dual $x^{*}$ of a point in $\mathbb{P}^{2}(\mathbb{F}_{p})$ is a line.
%$$(x^{*})^{p^{2}}=1 \iff x^{p^{2}-1}x=1\iff x^{(p-1)(p+1)}x=1 \iff x^{p+1}=1,$$ and
 There are $p+1$ points on every line and $p+1$ lines through every point.
 By induction on $n$, we easily see that the number of points in $\mathbb{P}^{n}(\mathbb{F}_{p})$ is $p^{n}+p^{n-1}+\ldots+p+1$.
 
 %and we count them over $\mathbb{F}_{q}$ providing a closed formula in terms of Stirling numbers.
 
% In the present note, we study genus $g$ projective plane algebraic curves that are $d$ sheeted coverings of $\mathbb{P}^{1}$ and we count them over $\mathbb{F}_{p}$. Each isomorphism class over $\mathbb{F}_{p}$ is counted with the reciprocal of the number of $\mathbb{F}_{p}-$automorphisms.
  
%In coordinates, if $\{b_{0},b_{1},b_{2}\}$ is a basis of $\mathbb{F}_{p^{3}}$ over $\mathbb{F}_{p}$, a line is represented by a linear equation: $x_{0}b_{0}+x_{1}b_{0}|\, x_{0},x_{1}\in \mathbb{F}_{p}$.

%\begin{defi}
Let $\mathbb{F}_{q}[x_{1},x_{2}, x_{3}]$ be the algebra of polynomials in 3 variables over $\mathbb{F}_{q}$.
Any irreducible homogeneous polynomial $f\in \mathbb{F}_{q}[x_{1},x_{2},x_{3}]$, defines the curve $\mathcal{C}_{f}$ on $\mathbb{P}^{2}(\mathbb{F}_{p})$ as the vanishing locus $V(f)$ of the polynomial $f$. %is nonsingular. 
More precisely $V(f)$ is the set of rational points of the curve $\mathcal{C}$ and $\langle f \rangle$ is the ideal generated by  $f$ in $\mathbb{F}_{q}[x_{0},x_{1},x_{2}]$. The degree $d$ of the curve is the degree of its polynomial, and assuming that the curve $C_{f}$ is smooth projective, the genus is $(d-1)(d-2)/2$. 
We define the number $N(q)$ of $\mathbb{F}_{q}-$rational points on the curve to be
$$N(q)=|\{(x_{0},x_{1},x_{2})\in \mathbb{P}^{2}(\mathbb{F}_{q})| f(x_{0},x_{1},x_{2})=0\}|.$$
It is a polynomial in $q$ with integer coefficients, whenever $q$ is a prime power.
The average number of $\mathbb{F}_{q}-$rational points %$\sharp \mathcal{C}_{f}(\mathbb{F}_{q})$ 
over all polynomials $f\in \mathbb{F}_{q}[x_{0},x_{1},x_{2}]$ defining smooth projective curves $\mathcal{C}_{f}\in \mathbb{P}^{2}(\mathbb{F}_{q})$ is $q+1$, (see \cite{BDFL}).

Let $M_{g,n}$ be the moduli space of curves of projective nonsingular curves of genus $g$ and $n$ marked points, and $\overline{M}_{g,n}$ its Deligne-Mumford compatification, corresponding to projective, connected, nodal curves of genus $g$ with $n$ marked point satisfying a stability condition. %Let $\overline{M}(\mathbb{F}_{q})$ denote the set of $\mathbb{F}_{q}$ points. 
For various ranges of $g$, $n$ and $q$, counting the number of points of $\overline{M}_{g,n}$ is feasible (see \cite{FP}).
%.% and the genus.

%\end{defi}

Consider points $p_{i}=[x_{i}:1]$ in $\mathbb{P}^{2}(\mathbb{F}_{q})$ and integers $m_{i}\geq 0$ such that $d|\, \sum_{i=1}^{n} m_{i}$, we study function field extensions of the rational function field $\mathbb{F}_{q}(x)$ given by an equation of the form:
 $$y^{n}=\prod_{i=1}^{d}(x-x_{i})^{m_{i}}.$$

We prove that  these field extensions correspond to plane algebraic curves $\mathcal{C}$ that are cyclic coverings of the projective line $\mathbb{P}^{1}(\mathbb{F}_{q})$. %where $f:=\prod_{i=1}^{n}(x-x_{i})^{m_{i}}.$
Let us denote by $F$ the function field of the curve $\mathcal{C}$ and by $\mathbb{P}(F)$ the set of places of $F$. For $P\in \mathbb{P}(F)$ we denote by $\mathcal{O}_{P}$ the valuation ring of $P$, and by $\nu_{P}$ the discrete valuation induced by $P$ in $F$. The set of places $\mathbb{P}_{1}(F)$ of $F$ of degree one are called rational places of $F$, i.e. $N(F)=|\mathbb{P}_{1}(F)|$.
In the present paper, we will concentrate in curves that are coverings of the projective line $\mathbb{P}^{1}(\mathbb{F}_{q})$ where $q$ is a power of a prime $p$. Moreover in Theorem \ref{Teo1}, we characterize all plane cyclic coverings of $\mathbb{P}^{1}(\mathbb{F}_{q})$ that are obtained by projecting its automorphism group to the known finite subgroups of $PGL(2,q)$ the automorphism group of the rational function field. Furthermore, we count them over $\mathbb{F}_{q}$. Each isomorphism class over $\mathbb{F}_{q}$ is counted with the reciprocal of the number of $\mathbb{F}_{q}-$automorphisms.

%We will assume that $\mathbb{F}_{q}$ is the full constant field of $F$
\section{Kummer extensions of the rational function field}
From now  $\mathbb{F}_{q}$ will be a field  with $q=p^{n}$ elements and $\mathcal{C}$ a non-singular, projective, irreducible curve defined over  $\mathbb{F}_{q}$, defined by the vanishing locus of a polynomial $F\in \mathbb{F}_{q}[x_{0},x_{1},x_{2}]$. Let $F/\mathbb{F}_{q}$ be the function field of $\mathcal{C}$. It is a finite algebraic field extension of $ \mathbb{F}_{q}(x)$ for some $x\in F$ that is transcendental over $\mathbb{F}_{q}$.
%We define the number $N(q)$ of $\mathbb{F}_{q}-$rational points on the curve to be
%$$N(q)=|\{(x_{0},x_{1},x_{2})\in \mathbb{P}^{2}(\mathbb{F}_{q})| F(x_{0},x_{1},x_{2})=0\}|.$$
%It is a polynomial in $q$ with integer coefficients, whenever $q$ is a prime power.

The number of points $N(q^{r})$%$\overline{\mathcal{C}}(\mathbb{F}_{q^{r}})$ 
on $\mathcal{C}$ over the extensions $\mathbb{F}_{q^{r}}$ of $\mathbb{F}_{q}$ is encoded in an exponential generating series, called the zeta function of $\overline{\mathcal{C}}$:
 $$Z(\mathcal{C},t)={\rm{exp}}\left( \sum_{r=1}^{\infty} N(q^{r})\frac{t^{r}}{r}\right).$$

%\begin{theorem}
The number of $\mathbb{F}_{q}$ rational points  of a nonsingular plane curve $\mathcal{C}$ of genus $g$ and degree $d$ satisfies the following bounds:
\begin{enumerate}
\item (Hasse-Weil bound) $N\leq q+1+2g\sqrt{q}$.
\item (St\"ohr-Voloch bound) $N\leq \frac{1}{2}(2g-2+(q+2)d).$
\item (Weil) $N(q^{2})\leq 1+q^{2}+2qg.$
\end{enumerate}
%\end{theorem}

The curve $\mathcal{C}$ is called $\mathbb{F}_{q^{2}}-$maximal if it attains the upper bound above; i.e. if one has
$$N(q^{2})=1+q^{2}+2q\cdot g.$$

The most well-known example of a $\mathbb{F}_{q^{2}}-$maximal curve is the so called hermitian curve $\mathcal{H}$, which can be given by the plane model:
$$x^{q+1}+y^{q+1}+z^{q+1}=0.$$
Thus many examples of $\mathbb{F}_{q^{2}}-$maximal curves arise by considering quotient curves $\mathcal{H}/G$, where $G$ is a subgroup of the automorphism group of $\mathcal{H}$.

%\begin{lemma} Let $Q(2,p)$ be the set of quadrics in $PG(2,p)$, that is, the varieties $V(g)$, where
%$g=a_{00}x_{0}^{2}+a_{11}x_{1}^{2}+a_{22}x_{2}^{2}+a_{01}x_{0}x_{1}+a_{02}x_{0}x_{2}+a_{12}x_{1}x_{2}$,  $a_{ij}\in \mathbb{F}_{p}$, then the cardinality:
%$$|Q(2,p)|=\frac{p^{6}-1}{p-1}$$
%\end{lemma}
%{\it Proof.} Observe that $(\mathbb{F}_{p^{2}})^{3}\cong \mathbb{F}_{p^{6}}$.

%Kummer extensions of the projective line over finite fields
%Any polynomial $\lambda(x,t)\in \mathbb{F}_{q}[x,t]$ of degree $n$ in $x$ defines a degree $n$ extension $\mathbb{F}_{q}[x]/(\lambda(x,t))$ of the rational function field $\mathbb{F}_{q}(t)$.

%The Galois group of the evaluated polynomial $f(x,\alpha)$ for any $\alpha \in \mathbb{F}_{q}$ is a subgroup of $G_{f}$. Every intermediate subfield $E$ of $F'/F$ is the function field of a smooth projective curve $X_{E}$ over $\mathbb{F}_{q}$, and if $E\subset E'$ are two such subfields, then there exists a $\mathbb{F}_{q}-$morphism $X_{E'}\rightarrow X_{E}$ of degree $[E':E]$. Intermediate subfields $E$ of $F'/F$, are in bijective correspondence with subgroups $S=Gal(F'/E)$ of $G_{f}=Gal(F'/F)$.
\begin{lemma}\label{lem1} Any finite Galois extension $F=\mathbb{F}_{q}(t)\hookrightarrow F'$ corresponds to a Galois covering $C\rightarrow \mathbb{P}^{1}(\mathbb{F}_{q})$. If the extension is cyclic then  $Gal(C/\mathbb{P}^{1}(\mathbb{F}_{q}))\cong \mathbb{Z}_{n}$, the cyclic group of order the degree $n$ of the minimal polynomial of the extension $F\hookrightarrow F'$. \end{lemma}

{\it Proof.} Let  $\mathbb{F}_{q}(t)$ be the rational function field of the projective line $\mathbb{P}^{1}(\mathbb{F}_{q})$. Any finite Galois extension $F=\mathbb{F}_{q}(t)\hookrightarrow F'$ corresponds to an inclusion of function fields $\mathbb{F}_{q}(t)\hookrightarrow \mathbb{F}_{q}(y,t)$, where $y$ satisfies a polynomial algebraic equation $f(x,t)$ over $\mathbb{F}_{q}[x,t]$, such that $f$ has degree $n>0$ in $x$, it is irreducible over $\mathbb{F}_{q}(t)$ and its zero locus defines an algebraic curve $C_{f}$ with function field $\mathbb{F}_{q}(y,t)$. The Galois group of the extension $G_{f}=Gal(F'/F)$ is a normal cyclic subgroup $\mathbb{Z}_{n}$ of the automorphism group ${\rm{Aut}}(C_{f})$ of the curve of order the degree $n$ of the field extension $F'/ F$. %minimal polyomial of the extension. % a quotient of the automorphism group ${\rm{Aut}}(C_{f})$ of the curve. %When $G_{f}$ contains a cyclic group $C_{m}$ such that the quotient curve $C_{f}/\mathbb{Z}_{m}$ has genus 0, then the curve is called a cyclic covering of $\mathbb{P}^{1}$.
%isomorphic to a quotient of the fundamental group of $\mathbb{P}^{1}$, thus it is a triangle group. 
Namely, the quotient curve $C_{f}/\mathbb{Z}_{n}$ has genus 0, and the curve is called a cyclic covering of $\mathbb{P}^{1}$.
\cqd

\begin{remark} 
The Galois group of the evaluated polynomial $f(x,\alpha)$ for any $\alpha \in \mathbb{F}_{q}$ is a subgroup of $G_{f}$. Every intermediate subfield $E$ of $F'/F$ is the function field of a smooth projective curve $X_{E}$ over $\mathbb{F}_{q}$, and if $E\subset E'$ are two such subfields, then there exists a $\mathbb{F}_{q}-$morphism $X_{E'}\rightarrow X_{E}$ of degree $[E':E]$. Intermediate subfields $E$ of $F'/F$, are in bijective correspondence with subgroups $S=Gal(F'/E)$ of $G_{f}=Gal(F'/F)$. %The projection to $t-$map defines a branched covering.
\end{remark}

\begin{example}
One can consider towers $\mathbb{F}_{q}(t^{\frac{1}{n}})$ or $\overline{\mathbb{F}}_{q}(t^{\frac{1}{n}})$, as $n$ varies through powers of the prime $p$ or through all integers not divisible by the characteristic of the ground field, that is $p$. The corresponding field extension $F=\mathbb{F}_{q}(t)\hookrightarrow F'=\mathbb{F}_{q}(t,y)$, where $y$ is a $n-$root of the polynomial $\sigma(u,t)=t^{n}-u \in F[u]$ and $u$ is transcendental over $\mathbb{F}_{q}$, is a finite cyclic Galois extension of degree $n$. Moreover $F'/F$ is an extension of Kummer type and $\mathbb{F}_{q}$ is the full constant field. %\cqd
\end{example}

\begin{remark} One of the main problems in coding theory is to obtain non-trivial lower bounds of the number $N(F_{i})$ of rational places of towers of function fields $\{F_{i}/\mathbb{F}_{q}\}_{i=1}^{\infty}$ such that $F_{i}\subsetneq F_{i+1}$.
\end{remark}

\begin{lemma}\label{L1} 
Any  Galois extension $\mathbb{F}_{q}(x)\hookrightarrow \mathbb{F}_{q}(x,y)$ of the rational function field, %after a birational transformation can be written in the form $
corresponding to a plane curve $C$ with affine model given by a polynomial equation of the form:
$$y^{n}=\prod_{i=1}^{s}(x-\rho_{i})^{d_{i}},\ \ \ d_{i}\in \mathbb{Z},$$ in some algebraically closed field containing $\mathbb{F}_{q}$, is  a cyclic plane covering of $\mathbb{P}^{1}(\mathbb{F}_{q})$, whose automorphism group contains a cyclic subgroup $\mathbb{Z}_{n}$. 
\end{lemma}

{\it Proof.} By Lemma \ref{lem1} any cyclic Galois cover of the projective line $\mathbb{P}^{1}(\mathbb{F}_{q})$ corresponds to a finite field extension $F\hookrightarrow F(y)$ with $y$  a $n-$root of the minimal polynomial $\sigma(t)=t^{n}-u \in F[u]$.  

Let $S=\mathbb{P}_{1}(F)\backslash \{ \rm{poles}\ \ \rm{and}\ \  \rm{zeroes}\ \ \rm{of}\ \ u \ \ in \ \ F\}.$ For any $P\in S$ we have that the polynomial $\overline{\sigma}_{P}(u,t):=t^{n}-u(P)$ factorizes in $\mathbb{F}_{q}[u,t]$ into pairwise distinct irreducible factors (see \cite{CT}), and by induction in the degree it follows the result.\cqd

\begin{remark}
\begin{enumerate}
\item The covering ramifies exactly at the places $x=\rho_{i}$, and the corresponding ramification indices are defined by $e_{i}=\frac{n}{(n,d_{i})}$, with $d_{i}$ the corresponding multiplicity of $\rho_{i}$. If $d=\sum_{i=1}^{s}d_{i}\equiv\, 0\, mod\, q$, then the place at $\infty$ does not ramify at the above extension. The only places of $F$ that ramify are the places $P_{i}$ that correspond to the points $x=\rho_{i}$.

\item If $(n,d_{1},\ldots, d_{s})=1$, then $F'$ is a Kummer extension of the rational function field $F$ of order $n$.

\item The genus of the function field $F'$ can be computed with the aid of the Riemann-Hurwitz formula to be  $g(F')=\frac{(n-1)(s-2)}{2}$.
\item The condition $(n,d_{i})=1$ is a stronger condition for all $i=1,\ldots, s$.

\end{enumerate}
\end{remark}

\begin{defi} A curve $\mathcal{C}$ is a $p-$cyclic cover of the projective line if and only if has a $g^{1}_{p}$ base point free linear system. If the linear system is unique, then the Galois cyclic group $Gal(\mathcal{C}/\mathbb{P}^{1})$ is normal in $G$. A sufficient condition for $g^{1}_{p}$ to be unique is the inequality: 
$$2\leq p \leq \frac{g}{2}+1.$$
\end{defi}

\begin{theorem} \label{Teo1} 
\begin{enumerate}
\item We can determine all cyclic Galois coverings of the projective line by considering every finite subgroup of the projective linear group $PGL(2,q)$, the automorphism group of $\mathbb{P}^{1}(\mathbb{F}_{q})$.
\item  Fix an odd prime number $p$, then the number of  $p-$sheeted coverings $y^{p}=f_{n}(x)$ of $\mathbb{P}^{1}(\mathbb{F}_{q})$ defined over $\mathbb{F}_{q}$ is given by the sum $\sum_{k=1}^{n}(q)_{k}$, where $(q)_{k}$ is the falling factorial polynomial $q\cdot (q-1)\ldots (q-(k-1))$, divided by the order of the affine transformation group of $\mathbb{A}^{1}=\mathbb{P}^{1}\backslash \infty$, that is $q^{2}-q$.

\end{enumerate}
\end{theorem}
{\it Proof.} %We observe that by Lemma \ref{L1}, 
Consider a cyclic Galois covering of the projective line $\mathbb{P}^{1}(\mathbb{F}_{q})$  given by a surjective morphism $\pi: \mathcal{Y}\rightarrow \mathbb{P}^{1}(\mathbb{F}_{q})$, where $\mathcal{Y}$ is a curve  with affine model defined by an equation $y^{p}=f(x)$ in  $\mathbb{F}_{q}[x,y]$, where $p$ is the degree of the corresponding field extension $\mathbb{F}_{q}(x)\hookrightarrow \mathbb{F}_{q}(x,y)$. %Then the automorphism group $Aut\,(\mathcal{Y})$ contains a cyclic subgroup $\mathbb{Z}_{p}$
Since the extension is cyclic, the invariant field $F^{0}:=(\mathbb{F}_{q}(x,y))^{\mathbb{Z}_{p}}$, by the action of the cyclic group, is the rational function field $\mathbb{F}_{q}(x)$. %and we call it $F^{0}$, 
 %Let us call $G$ to the automorphism group of the curve $\mathcal{Y}$. 
In particular, there is a map from the Galois group $G$ of the extension field $F'/F^{0}$, that projects into the subgroups of $PGL(2,q)$, the automorphism group of the rational function field
%such that $\mathcal{Y}/G$ is isomorphic to $\mathbb{P}^{1}(\mathbb{F}_{q})$, being $G$ a finite subgroup of $PGL(2,q)$ the automorphism group of the rational function field
 that acts on $\mathbb{P}^{1}(\mathbb{F}_{q})$ by:
$$\left(\begin{array}{ll} a &  b \\
c & d \end{array}\right) : z \rightarrow \frac{az+b}{cz+d}, \ \ \  A\in PGL(2,q).$$ %We call it
%$X_{\tau}$,
This corresponds to the fact that any cyclic Galois covering $\pi: \mathcal{Y}\rightarrow \mathbb{P}^{1}(\mathbb{F}_{q})$ arises as the quotient curve of an algebraic projective curve defined over $\mathbb{F}_{q}$ by a finite subgroup of $PGL(2,q)$, that it turns out to be the fundamental group of the covering curve $\mathcal{Y}$.

%the Galois group of the field extension $F'/F^{0}$ is isomorphic to the fundamental group of $\mathbb{P}^{1}$.

If $G=PSL_{2}(q)$ (respectively, $PGL_{2}(q)$) with $q\leq 5$ (respectively, $q\leq 4$), then $G$ is isomorphic to one of $S_{3}, S_{4}, A_{4}$ or $A_{5}$.
If $G=PSL_{2}(q)$ with $5<q=p^{e}$, for some prime number $p$ and some positive integer $e$, $G$ is a quotient with torsion free kernel of a certain triangle group:
$$T_{r,s,t}=<x,y,z: \ \ x^{r}=y^{s}=z^{t}=xyz=1>,$$
where $r,s,t$ are integer numbers such that $\frac{1}{r}+\frac{1}{s}+\frac{1}{t}<1$, the group $T_{r,s,t}$ is called a hyperbolic triangle group. %If $\frac{1}{r}+\frac{1}{s}+\frac{1}{t}>1$ then $T_{r,s,t}$ is a finite group and it is either dihedral or isomorphic to one of $A_{4}, A_{5}$ or $S_{4}$, (see \cite{SG}).

In order to prove the second part ot the Theorem, we need to count all the polynomials $f_{n}(x)$ that define a $p-$sheeted covering of the projective line, by counting polynomials depending on the number of different roots. We implicitly assume, that our polynomial $f_{n}(x)$ decomposes into linear factors. Otherwise we work over $\overline{\mathbb{F}}_{q}[x]$, where $\overline{\mathbb{F}}_{q}$ denotes the algebraic closure of the finite field $\mathbb{F}_{q}$.

Thus the number of separable polinomials of degree $n$ over $\mathbb{F}_{q}$ is $q^{n}$. The number of monic polynomials with $n-1$ different roots is 

\noindent $q(q-1)(q-2)\ldots (q-n+1)$, that is known as the falling factorial polynomial $(q)_{n}$ and it is the generating function for the signed Stirling numbers: \footnote{The Stirling number of the first kind $\left[\begin{matrix} n \\ k \\ \end{matrix} \right ]$  is the number of ways to partition a set of $n$ objects into $k$ non-empty subsets and the signed Stirling number  of the first kind denoted by $s(n,k)$ is recovered from the unsigned one, by the rule $s(n,k)=(-1)^{n-k}\left[\begin{matrix} n \\ k \\ \end{matrix} \right ]$.}

$q_{(n)}=q(q-1)(q-2)\ldots (q-n+1)=\sum_{k=0}^{n}s(n,k)q^{k}$.
%, the combinatorial number ${p\choose n-1}$.

One can immediately see, that $(q)_{1}=q$, because the number of polynomials with one different root is $q$, independently of  the degree $n$ of the polynomial $f(x)$. The number of polynomials with $k$ different roots, where $1\leq k \leq n$, is given by the falling factorial polynomial $(q)_{(k)}$.

Since only $\infty\in \mathbb{P}^{1}$ is distinguished, we must further divide by the group of affine transformations of $\mathbb{A}^{1}=\mathbb{P}^{1}\backslash \infty$. Since the order of the transformation group is $q^{2}-q$, we have that the number of $q$ curves over $\mathbb{F}_{q}$ expressible as $d-$sheeted coverings of $\mathbb{P}^{1}$ is the falling factorial polynomial $(q-2)_{n}$. %what is the role of $d$ here

 \cqd

\begin{remark} The case in which $\mathbb{P}^{1}(\mathbb{F}_{q})$ is covered by a curve $\mathcal{C}$ determined by a finite subgroup $G<PGL(2,q)$ such that the corresponding projection to the quotient surface $\mathcal{C}\rightarrow \mathcal{C}/G\cong \mathbb{P}^{1}(\mathbb{F}_{q})$ is branched over 3 points plays an important role in the study of Beauville surfaces, (see \cite{SG}).
\end{remark}

\begin{remark} Theorem \ref{Teo1} gives a motivation for the study of quotient curves arising from subgroups of $PSL(2,q)$ or $SL(2,q)$. In particular, given a polynomial $f\in \mathbb{F}_{q}[x]$ whose degree $d$ is not divisible by $p$, it applies not only to cyclic covers of $\mathbb{P}^{1}(\mathbb{F}_{q})$ defined by $y^{p}=f(x)$ but also to Artin-Schreier covers $C_{f}$ of $\mathbb{P}^{1}$ with affine model $y^{p}-y=f(x)$.
\end{remark}

\subsubsection*{Acknowledgments}
We thank professor Michael Zieve for comments and suggestions during the preparation of the paper.% and the organizers of the conference %Cuartas Jornadas en Teor\'ia de N\'umeros for gives us the opportunity to present our work.

\end{document}